\documentclass[11pt]{article}
\usepackage[utf8]{inputenc}
\usepackage[english]{babel}
\usepackage{amsmath,amssymb,amsthm,amsfonts,mathrsfs}
\usepackage{graphicx}
\usepackage{tikz}
\usepackage{xcolor}
\usepackage{color}
\usepackage{xspace}
\usepackage{calrsfs}
\usepackage{algorithm}
\usepackage{algpseudocode}
\usepackage{adjustbox}
\usepackage{xspace}

\theoremstyle{definition}

\theoremstyle{plain}
\newtheorem{thm}{Theorem}[section]

\newtheorem{example}[thm]{example}
\newcommand{\todocmd}[1]{\small{\textcolor{red}{#1}}}

\newcommand{\todo}[1]{\ifthenelse {\boolean{showComments}} {\todocmd{#1}} {}}

\newcommand{\Sys}{R}
\newcommand{\Syshat}{\hat R}

\newcommand{\dd}[1]{\frac{\partial}{\partial #1}}

\newcommand{\BK}{\sc{BK}}

\newcommand{\DE}{DE\xspace}
\newcommand{\DEs}{DEs\xspace}
\newcommand{\PDE}{PDE\xspace}
\newcommand{\ODE}{ODE\xspace}

\newcommand{\DPS}{DPS\xspace}

\newcommand{\shatn}{\hat{s}}
\newcommand{\thatn}{\hat{t}}
\newcommand{\vhatn}{\hat{v}}
\newcommand{\xh}{\hat{x}}

\newcommand{\uh}{\hat{u}}
\newcommand{\thn}{\hat{t}}

\newcommand{\alg}[1]{{\sc #1}\xspace}
\newcommand{\maple}[1]{{\tt #1}\xspace}
\newcommand{\inputs}[1]{{\it #1}\xspace}
\newcommand{\MapDE}{\maple{MapDE}}
\newcommand{\Source}{\inputs{Source}}

\newcommand{\Target}{\inputs{Target}}
\newcommand{\TargetClass}{\inputs{TargetClass}}
\newcommand{\Map}{\inputs{Map}}

\newcommand{\Maple}{\maple{Maple}}
\newcommand{\LAVF}{\maple{LAVF}}

\newcommand{\rif}{\alg{rif}}

\newcount\dotcnt\newdimen\deltay
\def\Ddot#1#2(#3,#4,#5,#6){\deltay=#6\setbox1=\hbox to0pt{\smash{\dotcnt=1
			\kern#3\loop\raise\dotcnt\deltay\hbox
			to0pt{\hss#2}\kern#5\ifnum\dotcnt<#1 \advance\dotcnt
			1\repeat}\hss}\setbox2=\vtop{\box1}\ht2=#4\box2}

\title{\large \bf Introduction of the MapDE algorithm for determination of \\ mappings relating differential equations}

\author{\bf\small\sc Zahra. Mohammadi$^{\dag}$, Gregory J. Reid \footnote{Corresponding author. Email:  reid@uwo.ca; zmohamm5@uwo.ca; Tracy.Huang@data61.csiro.au } and S.-L. Tracy Huang$^{\ddag}$ \\ \small\em{$^{\dag, *}\,$Department of Applied Mathematics, University of Western Ontario, Canada} \\ \small\em{$^{\ddag}\,$Data61, CSIRO, Canberra ACT 2601, Australia} }

\date{2019}

\begin{document}	
\maketitle

\begin{abstract}
	This paper is the first of a series in which we  develop exact and approximate algorithms for mappings of systems of differential equations.  Here we introduce the \MapDE algorithm and its implementation in Maple, for mappings relating differential equations.  We consider the problem of how to algorithmically characterize, and then to compute mappings of less tractable (\Source) systems $\Sys$ to more tractable (\Target) systems $\Syshat$ by exploiting the Lie algebra of vector fields leaving $\Sys$ invariant. 
	Suppose that $\Sys$ is a (\Source) system of (partial or ordinary) differential equations with independent variables $x = (x^1, x^2, \ldots , x^n ) \in \mathbb{C}^n$ and dependent variables $u = (u^1, \ldots , u^m)\in \mathbb{C}^m$.  Similarly suppose $\Syshat$ is a (\Target) system in the variables $(\hat{x}, \hat{u}) \in \mathbb{C}^{n+m}$.  For systems of exact differential polynomials $\Sys$, $\Syshat$  our algorithm \MapDE can decide, under certain assumptions, if there exists a local invertible mapping  $ \Psi(x,u) = (\hat{x}, \hat{u} ) $ that maps the \Source system $\Sys$ to the \Target $\Syshat$.
	We use a result of Bluman and Kumei who have shown that the mapping $\Psi$ satisfies infinitesimal (linearized) mapping equations that map the infinitesimals of the Lie invariance algebra for $\Sys$ to those for $\Syshat$.  
	
	\MapDE involves applying the differential-elimination algorithm to the defining systems for  infinitesimal symmetries of $\Sys$, $\Syshat$, and also to the nonlinear mapping equations (including the Bluman-Kumei mapping subsystem); returning them in a form which includes 
	its integrability conditions and for which an existence uniqueness theorem is available.
	Once existence is established, a second stage can determine features of the map, and some times by integration, explicit forms of the mapping.  Examples are given to illustrate the algorithm.
	
	Algorithm \MapDE also allows users to enter broad target classes instead of a specific system $\Syshat$.
	For example we give an algorithmic approach that avoids the integrations of the Bluman-Kumei approach
	where \MapDE can 
	determine if a linear differential equation $\Sys$ can be mapped to a linear constant coefficient differential equation.
	
\end{abstract}

\vspace{0.10in}

\vspace{0.1in}
\noindent {\it Categories and Subject Descriptors}: I.1.2; I.1.4.

\section{Introduction}

\label{sec:Intro}

This paper is the first of a series in which we explore algorithmic aspects of mappings of
differential equation systems that transform differential equations (\DEs) to \DEs.  Naturally this exploration includes symmetry transformations 
-- transformations of a \DE to itself, and also equivalence transformations where one member of a class of \DEs is mapped 
to another member of the class.  
In this paper we introduce the algorithm \MapDE for characterizing mappings between \DEs.  For algorithmic implementation we restrict our treatment to differential polynomial systems ({\DPS}), systems which are polynomially nonlinear functions of their derivatives and dependent variables; with coefficients from 
some computable field (e.g. $\mathbb{Q}$).

In earlier work we developed approximate methods for determination of approximate Lie symmetry algebra of \DEs  \cite{BLRSZ:2004, LisHG:Sym}.  
A key motivation for our current work, is how to practically use such approximate methods.   We see determination of approximate mappings of {\DPS}, to be explored later in
this series, as a practical way in which to exploit such approximate symmetry information.  Our interest in mappings was also motivated by recent work \cite{LGM101:LG}, which used Reid \cite{Rei91:Fin} on the algorithmic determination of structure of Lie algebras of symmetries of {\DE}, to give an algorithm to determine the existence of mappings exactly linearizable {\ODE}.
We give an algorithmic implementation of the methods of Bluman and Kumei \cite{Blu10:App, BluKu109:DiEq} for exploiting the Lie symmetries of a system in the determination of mappings between \DEs.

In particular in this paper we introduce an algorithm for such mappings in the presence of symmetry.  The algorithm {\MapDE} is implemented as part of Huang and Lisle's \LAVF object-oriented Maple package~\cite{LisH:Alg}.
We give examples to illustrate the algorithm and compare it with the approach of Bluman and Kumei.
We extend the algorithm, to determine the existence of a mapping from linear DE, to linear constant coefficient \DE,
avoiding the heuristic integrations of Bluman and Kumei's approach.


We consider systems of (partial or ordinary) differential equations with $n$ independent variables and 
$m$ dependent variables. Suppose $\Sys$ has independent variables $x = (x^1, \ldots ,x^n )$ and dependent variables $u = (u^1, \ldots , u^m)$ and 
$\Syshat$ has independent variables $\hat{x} = (\hat{x}^1, \ldots , \hat{x}^n )$ and dependent variables $\hat{u} = (\hat{u}^1, \ldots , \hat{u}^m)$.
In particular we consider local analytic mappings $\Psi$:
$(\hat{x}, \hat{u} ) = \Psi (x,u) = (\psi (x,u), \phi(x,u) )$, so that $\Sys$ is locally and invertibly mapped to $\Syshat$:
\begin{equation}
\label{eq:MapTrans}
\hat{x}^j  =   \psi^j (x,u), \qquad 
\hat{u}^k =   \phi^k (x,u) 
\end{equation}
where $j = 1,\ldots, n$ and $k = 1,\ldots, m$.
The mapping is locally invertible so the determinant of the Jacobian of the mapping is nonzero: 
\begin{equation}
\label{eq:Jac}
\mbox{Det Jac}(\Psi) = \mbox{Det} \frac{\partial (\psi, \phi )}{\partial (x, u)} \not = 0,
\end{equation}
where $\frac{\partial (\psi, \phi )}{\partial (x, u)} $ is the usual Jacobian $(n+m) \times (n+m)$ matrix of first order derivatives of the $(n+m)$ functions $(\psi, \phi)$ with respect to the $(n+m)$ variables $(x, u)$.
Note throughout this paper, we will call $\Syshat$ the 
\Target system of the mapping, which will generally have some more desirable features than $\Sys$, which we call the \Source system.

Algorithms for existence of such mappings, and methods for their explicit construction, is the topic of this paper. 
A very general approach to such problems, Cartan's famous Method of Equivalence \cite{PetOlv107:Sym},
finds invariants, that label the classes of systems, equivalent under the pseudogroup of such mappings.
The fundamental importance of such equivalence questions, and the associated demanding computations 
has attracted attention from symbolic computation researchers.
For example, Neut, Petitot and Dridi \cite{NeutPetitotDridi2009}, implemented Cartan's method for {\ODE} and certain classes of {\PDE} of finite type (i.e.\ with finite dimensional solution space).
Olver and collaborators developed a new version of Cartan's moving frames \cite{Fel99:Mov}.
Valiquette \cite{Valiquette13:LocEquiv} applied this method to equivalence problems and 
further results are given by Arnaldson \cite{Arnaldsson17:InvolMovingFrames}.
Also see \cite{And12:New} which introduces the \maple{DifferentialGeometry} package, available in {\Maple}
and has been applied to equivalence problems. 
Also see 
\cite{Hub09:Dif, Man10:Pra} for approaches to the non-commutative calculus that results in calculations.
Underlying these calculations, is that overdetermined \PDE systems, with some 
non-linearity, are required to be reduced to forms that enable the statement of a local existence and uniqueness
theorem (such include passive and involutive forms).  See  \cite{GolKonOvcSza09:complexityDiffElim} for estimates of complexity of such methods, which indicate their difficulty.


Our initial approach, is fairly direct, and exploits the linearity of the Bluman-Kumei mapping equations.
It also is motivated by in the longer term, we wish to include invariant differential operators and using the newly
developed methods of Numerical Jet Geometry, to investigate approximate equivalence.

In contrast, Bluman and Kumei \cite{BluKu109:DiEq} consider a narrower class of mapping problems, which is focused on the case where the \Target system, is uniquely characterized in terms of its Lie symmetry invariance algebra (See \cite{Blu10:App, BluKu109:DiEq}).
In this article we will implement an algorithm based on the Bluman and Kumei approach.

Suppose that a \Source system, has an associated Lie symmetry algebra, together with its defining system.
Such infinitesimal Lie {\em point} symmetries for $\Sys$ are found by seeking
vector fields
\begin{equation}
\label{eq:symmOp}
V = \sum_{i=1}^n \xi^i(x,u) \dd{x^i} + \sum_{j=1}^m \eta^j(x,u)\dd{u^j}
\end{equation}
whose associated one-parameter group of transformations
\begin{align}
\label{eq:LieInfTrans}
{x^*}   &= x +  \xi(x,u) \epsilon + O(\epsilon^2) \nonumber\\
{u^*}   &= u + \eta(x,u) \epsilon + O(\epsilon^2)
\end{align}
which away from exceptional points preserves the jet locus of such systems - mapping solutions to solutions.
See \cite{BluKu111:Sym, BluKu112:Sym} for applications.
The {\em infinitesimals} $(\xi^i,\eta^j)$ of a symmetry vector field~\eqref{eq:symmOp} for a system of \DEs are found by solving an associated system of linear homogeneous defining equations (or determining equations) for the infinitesimals.
The defining system is derived by an explicit algorithm, for which numerous computer
implementations are available \cite{Car00:Sym, Che07:GeM, Roc11:SAD}.
Similarly we suppose that the \Target admits symmetry vector fields

\begin{equation}
\label{eq:vf}
\hat{V} = \sum_{i=1}^n \hat{\xi}^i(\hat{x},\hat{u}) \dd{\hat{x}^i} + \sum_{j=1}^m \hat{\eta}^j(\hat{x},\hat{u})\dd{\hat{u}^j}
\end{equation}
in the \Target infinitesimals $(\hat{\xi}, \hat{\eta})$.
Computations with defining systems of both systems will be essential in our approach.
We have implemented our algorithms in Huang and Lisle's powerful object oriented {\LAVF}, {\Maple} package~\cite{LisH:Alg}.

In \S \ref{sec:PreMapEqs} we give preliminaries and the Bluman-Kumei Mapping equations, together with a simple illustrative example.
In \S \ref{sec:Algorithms} we describe our core algorithm \MapDE, which takes $\Sys$ and $\Syshat$ as input, and 
returns the reduced involutive form \rif-form for $\Psi$, establishing existence, non-existence of the mapping.
Once existence of the mapping is established, a further phase, is to try to obtain an explicit form for the mapping by integrating the mapping equations.  We treat two cases: one in which $\Sys$ and $\Syshat$ are specified and another where
\maple{TargetClass = ConstantCoeffDE}.
In \S \ref{sec:Examples} we give examples of application \MapDE and conclude with a discussion in \S \ref{sec:Discussion}.
\section{Preliminaries \& Mapping Equations}

\label{sec:PreMapEqs}

For an algorithmic treatment, we limit the systems considered to being differential polynomials, with coefficients from 
computable subfield of $\mathbb{C}$ (e.g.\ $\mathbb{Q}$).  Some non-polynomial systems can be converted to differential 
polynomial form by the use of the Maple command, \maple{dpolyform}.

In the geometric approach to \DEs centers around the jet locus, where the derivatives are regarded as formal variables
and a map to polynomials, in our case, where the tools of algebraic geometry can be used.  In general systems of polynomial
equations and inequations must be considered (differences of varieties).
The union of prolonged graphs of local solutions is a subset of the jet locus in $J^q$, the jet space of order $q$.
For details concerning the Jet geometry of \DEs see \cite{Olv93:App,Sei10:Inv}. 

\begin{example}
	\label{ex:SimpleEx}
	Consider the famous Black Schole's equation which is fundamental in financial applications \cite{MasMa119:ExBlacSch}, we will use as an introductory simple example:
	\begin{equation}
	\label{eq:BS}
	\frac{\partial}{\partial t} v + \frac{s^2}{2} \left(\frac{\partial}{\partial s} \right)^2 v + s \frac{\partial}{\partial s} v - v = 0
	\end{equation}
	By inspection this equation has the obvious symmetry of translation in $t: t^* = t + \epsilon$ and scaling in 
	$s:  s^* = b s$.  Moreover the infinitesimal form of these symmetries (\ref{eq:LieInfTrans}) is generated 
	by the operators (\ref{eq:symmOp}) given by $\frac{\partial}{\partial t}$ and $s\frac{\partial}{\partial s}$.
	Since these vector fields obviously commute, it is natural to map to new coordinates in which:
	\begin{equation}
	\label{eq:BSops}
	\frac{\partial}{\partial t} = \frac{\partial}{\partial \hat{t}}, \; \; \; s\frac{\partial}{\partial s}  = \frac{\partial}{\partial \hat{s}}
	\end{equation}
	So by trivial integration the transformation $\hat{s} = \int \frac{ds}{s} = \log(s)  + c_1 $,  $\hat{t} = t + c_2$ should map the Black Schole's Equation into an equation invariant under two commuting translations, i.e. to a constant coefficient equation.  Indeed by inspection we find:
	\begin{equation}
	\label{eq:BShat}
	\frac{\partial}{\partial \hat{t}} v + \frac{1}{2} \left(\frac{\partial}{\partial \hat{s}} \right)^2 v - v = 0
	\end{equation}
	which is the famous Black-Schole's transformation of \eqref{eq:BS} to the backwards heat equation.
	This example simply illustrates that there can be a strong connection between symmetries admitted by an equation
	and mappings of the equation to convenient forms.
\end{example}
Indeed this illustrates the key idea of Bluman-Kumei's method for determining when a linear differential equation (\DE) in $n$ independent variables can be mapped to a (\Target) linear constant coefficient \DE:  that the \Target admits 
$n$ commuting translations.  Geometrically the \Source must correspondingly admit a subalgebra of its Lie symmetry 
algebra consisting of $n$ commuting symmetries (that act transitively on the space of its independent variables).

One can try to devise an algorithm for determining such symmetries explicitly.  In general this involves integrating systems of overdetermined \PDE, and, though advantageous in many applications, no general algorithm is known for this task.  In our 
paper we describe algorithms using a finite number of differentiations and eliminations, and no integrations, that guarantees the 
algorithmic determination of the existence of such transformations.
The computer algebra system {\Maple} has several excellent such differential elimination algorithms, and also excellent algorithms for generating the 
linearized equations for symmetries.

\begin{example}
	For the Black-Schole's Equation~\eqref{eq:BS}, the defining system for the infinitesimal symmetry operator has form
	$ \sigma(s,t,v)\dd{s} + \tau(s,t,v) \dd{t} + \eta(s,t,v) \dd{v} $.  Here comparing with \eqref{eq:vf} yields $(\xi^1, \xi^2)=(\sigma, \tau)$ and
	$\eta^1 = \eta$. The automatically generated unsimplified 
	system of defining equations for infinitesimal symmetries is:
	\begin{eqnarray}
	&& \tau_{{s}}=0,\tau_{{v}}=0,\tau_{{v,v}}=0, \sigma_{{v,v}}=0,\tau_{{s,v}}{s}^{2}-2\, \sigma_{{v}}=0,   \nonumber  \\
	&&\eta_{{v,v}}s-2\, \sigma_{{s,v}}s-4\, \sigma_{{v}}=0,  \nonumber   \\
	&& -2\,\eta_{{s,v}}{s}^{2}+ \sigma_{{s,s}}{s}^{2}+2\, \sigma_{{s}}s+6\, \sigma_{{v}}v-2\, \sigma+2\, \sigma_{{t}}=0,   \\
	&&\tau_{{s,s}}{s}^{3}-2\,\tau_{{s}}{s}^{2}+2\,\tau_{{v}}sv+2\,\tau_{{t}}s-4\, \sigma_{{s}}s+4\, \sigma=0,   \nonumber  \\
	&&-\eta_{{s,s}}{s}^{3}+2\,\eta_{{s}}{s}^{2}-2\,\eta_{{v}}sv+4\, \sigma_{{s}}sv+2\,\eta\,s-2\,\eta_{{t}}s-4\,v \sigma=0   \nonumber 
	\end{eqnarray}
	
	Application of a differential-elimination algorithm to this system augmented with $\eta_v = \eta/v$ yields:
	\begin{eqnarray}
	\label{eq:rifBS}
	&& \eta_{{s}}={\frac {v \left( 3\,\tau_{{t}}s+4\, \sigma_{{t}} \right) 
		}{4{s}^{2}}},\eta_{{t}}={\frac {v \left( 17\,\tau_{{t}}s-2\,\tau_{
				{t,t}}s+12\, \sigma_{{t}} \right) }{8s}},\eta_{{v}}={\frac {\eta}{v}},   \nonumber \\
	&& \tau_{{t,t,t}}=0, \sigma_{{t,t}}=0, \sigma_{{s}}={\frac {\tau_{{t}}s+2\, \sigma}{2s
	}},\tau_{{s}}=0, \sigma_{{v}}=0,\tau_{{v}}=0
	\end{eqnarray}
	Application of \rif's initial data algorithm  yields:
	\begin{eqnarray}
	\label{eq:IDBS}
	&&\eta(s_0,t_0,v_0)=c_1, \quad \tau(s_0,t_0,v_0) =  c_2,\quad \tau_t(s_0,t_0,v_0) =  c_3,  \nonumber \\
	&&\tau_{tt} (s_0,t_0,v_0)= c_4,\quad  \sigma(s_0,t_0,v_0) =  c_5,\quad \sigma_{t}(s_0,t_0,v_0) =  c_6 \qquad
	\end{eqnarray}
	The key aspect relevant for our paper is that (\ref{eq:rifBS}) and (\ref{eq:IDBS}) are obtained with algorithmic operations and
	in particular without integration.  In addition further algorithms from the {\LAVF} package can determine the structure of
	its Lie Algebra.  Indeed we find a two dimensional abelian subalgebra from that output, a necessary condition for the existence of a map of the Black-Schole's equation to a constant coefficient equation.
	
\end{example}

\subsection{Mapping Equations}
Assuming existence of a local analytic invertible map $\Psi = (\psi, \phi)$ between the \Source system $\Sys$ and the \Target system $\Syshat$ and applying it to the infinitesimals $(\hat{\xi}, \hat{\eta})$ yields what we will call the \textit{Bluman-Kumei (BK) mapping equations}:
\begin{align}
	\label{eq:BK}
	\sum_{i=1}^n \xi^i(x,u)  \frac{\partial \psi^k}{\partial x^i} + \sum_{j=1}^m \eta^j(x,u)\frac{\partial \psi^k}{\partial u^j}
	&= \hat{\xi}^k (\hat{x}, \hat{u})  \nonumber   \\
	\sum_{i=1}^n \xi^i(x,u) \frac{\partial \phi^\ell}{\partial x^i}+\sum_{j=1}^m \eta^j(x,u)\frac{\partial \phi^\ell}{\partial u^j} 
	&=	 \hat{\eta}^\ell (\hat{x}, \hat{u}) 
\end{align}
where $1 \leq k \leq n$ and $1 \leq \ell \leq m$.
See Bluman and Kumei \cite{BK105:BluKu,Blu10:App} for details and generalizations (e.g. to contact transformations). Note that all quantities on the LHS of the BK mapping equations~\eqref{eq:BK} are functions of $(x,u)$ including $\phi$ and $\psi$.

\begin{example}
	\label{ex:SimpleExMapEqs}
	We informally illustrate the {\BK} mapping equations on Example \ref{ex:SimpleEx}.
	Here we follow an approach based on heuristic integration of the symmetry defining system.
	Indeed the defining system is easily integrated to find the full $6$ dimensional Lie symmetry algebra.
	And among the basis of symmetries the reader can easily find the two operators previously by inspection:
	\begin{equation}
	\label{eq:BSabelianOps}
	\frac{\partial}{\partial t}, \quad s\frac{\partial}{\partial s}
	\end{equation}
	which implies that mapping has form: $\hat{s} = \psi^1 (s,t)$, $\hat{t} = \psi^2 (s,t)$, $\hat{v} = \phi(s,t) = v$.
	When the corresponding coefficients of (\ref{eq:BSabelianOps}) are substituted into the BK system (\ref{eq:BK}) we get:
	\begin{eqnarray}
	s \dd{s} \psi^1 (s,t) = 1,\quad s \dd{s} \psi^2 (s,t) = 0 \nonumber  \\
	\dd{t} \psi^1 (s,t) = 0, \quad    \dd{t} \psi^2 (s,t) = 1 
	\end{eqnarray}
	which yields by simple integration the same result as before for the mapping of the Black-Schole's to constant coefficient:
	\begin{equation}
	\label{eq:BSmap}
	\hat{s} = \psi^1 (s,t) = \log(s), + c_1, \quad \hat{t} = \psi^2 (s,t) = t + c_2
	\end{equation}
	Indeed this integrating and breaking down into a basis, is the method used by Bluman and Kumei.
	However it does not yield an algorithm, since it depends on heuristic integration.
	
	Finally we mention, that we are not opposed to integration, and in fact, a combination of integration and the algorithmic 
	methods of this article, are probably a preferable way to proceed in practice.

\end{example}


Let $S$, $\hat{S}$ denote the symmetry defining systems for the \Source system $\Sys$ and the \Target system $\Syshat$ respectively, with corresponding Lie symmetry algebras $\mathcal{L}$ and $\hat{\mathcal{L}}$.
If an invertible map $\Psi$ exists mapping $\Sys$ to $\Syshat$ then it most generally 
depends on $\dim(\mathcal{L}) = \dim(\hat{\mathcal{L}})$ parameters.
But we only need one such $\Psi$.  So reducing the number of such parameters, e.g.\ by restricting to 
a Lie subalgebra $\mathcal{L}^{'}$ of $\mathcal{L}$ with corresponding Lie subalgebra $\hat{\mathcal{L}}^{'}$ of $\hat{\mathcal{L}}$
that still enables the existence of such $\Psi$, is important in reducing the computational
difficulty of such methods.  
We will use the notation $S',\hat{S}'$ denote the symmetry defining systems of  Lie sub-algebras $\mathcal{L}^{'} $,
$\hat{\mathcal{L}}^{'}$ respectively.
See \cite{Blu10:App,PetOlv107:Sym} 
discussion on this matter.

\begin{example}
	\label{ex:restrict2subalg}
	The mapping of the linear Black-Schole's equation~\eqref{eq:BS} to a constant coefficient linear equation; we exploited the existence of a two dimensional abelian subalgebra.  Indeed if we are lucky enough to identify this subalgebra immediately, then it gives very simple mapping equations with only two parameters. In the general algorithm for mapping linear equations to constant coefficient equations we described later, we can first bring such
	equations to homogeneous form.  Restricting to symmetries, and mappings that retain the homogeneous form, can
	be imposed by restricting to symmetries with $\{ \eta_v = \eta/v, \sigma_v = 0, \tau_v = 0 \}$, which has a finite $6$ parameter Lie group of symmetries.
	Thus we rejection the unhelpful infinite super-position subgroup
	as we did in the Black-Schole's example earlier.  For more details see Bluman et al.~\cite{Blu10:App}.
	
\end{example}

\subsection{Algorithms EquivDetSys and DimEquivTest}

\label{sec:EquivDetSysDimEquivTest}
\noindent
With the \Source system $\Sys$, the \Target system $\Syshat$ and the mapping $\Psi$, the algorithm $\maple{MapEqs}(\Sys, \Syshat, \Psi)$ is to return the full non-nonlinear defining equations for mappings from $\Sys$ to $\Syshat$ which are invertible (i.e.\ $\mbox{Det Jac}(\Psi) \not = 0$).  As preparation for the description of this algorithm we introduce the following algorithm.

\noindent
\textbf{EquivDetSys}$(\Sys, \Syshat)$:  This is {\Maple} implementation of returning the nonlinear {\DPS} for invertible mappings $\Psi$ from $\Sys$ to $\Syshat$.  Our implementation currently requires that $\Sys$, $\Syshat$ are in solved form
for their leading derivatives with respect to a ranking graded by total differential order; though this could be weakened in
the future.  Then {\Maple}'s general purpose routine for changing variables $\maple{dchange}$ is applied, yielding expressions
in the parametric derivatives of $\Sys$.  Setting coefficients of independent powers of the parametric derivatives to zero, 
together with $\mbox{Det Jac}(\Psi) \not = 0$ simplified with respect to $\Psi$ yields the nonlinear determining system
for $\Psi$.  This construction is well-known (indeed it is used in \cite{LGM101:LG} in the special case of mappings linearizing {\ODE}).  
However the nonlinear overdetermined systems are challenging to compute due to the expansion of
determinants as the number of variables and differential order of $\Sys, \Syshat$ increase.

Our approach in this paper, is to take advantage of such linearized infinitesimal information, available from Lie symmetries and in particular via the 
{\BK} equations, which are linear in the mapping variables.   Then if necessary, at the end of {\MapDE} apply
\textbf{EquivDetSys}, which can be much simplified by the earlier computed conditions in ${\Psi}$.

Also we employ a number of efficient preliminary tests that can some times quickly determine if $\Sys$ and $\Syshat$ 
are not equivalent via $\Psi$.

\noindent
\textbf{DimEquivTest}$(\Sys, \Syshat)$:  Differential-elimination algorithms such as those in the packages
$\maple{RIF}$, \maple{DifferentialAlgebra} and \\
\maple{DifferentialThomas} 
allow a determination of a coordinate dependent description of initial data, and using that the determination of the 
coordinate independent quantities $\dim(\Sys)$, $\dim(\Syshat)$.  Thus a quick first test applied 
by \maple{DimEquivTest} is $\dim(\Sys) = \dim(\Syshat)$.

If the input ranking is graded first by 
total derivative order, then further dimension invariants can be derived from that initial data:
which are the number of parametric derivatives at each derivative order $n$ (determining the Differential Hilbert Series).
\maple{DimEquivTest} tests the equality of these invariants up to the maximum involutivity order for $\Sys, \Syshat$.
One further invariant is the number of arbitrary functions of the maximum number of independent variables appearing 
in the initial data.  For background information see \cite{Sei10:Inv}.
\section{MapDE Algorithm}
\label{sec:Algorithms}
In this section we describe algorithms for mapping a system $\Sys$ to $\Syshat$.  

In \S\ref{sec:MapDE} we describe \MapDE for a specific \Source system $\Sys$ and 
specific \Target system $\Syshat$. 
In \S\ref{sec:ToConstantCoeffDE} we give a description of {\MapDE} for a linear input equation and
a class of \Target systems (where the \Target is constant coefficient linear equation).

\subsection{The MapDE Agorithm for specific $\Sys$ and $\Syshat$}

\label{sec:MapDE}

First we describe \MapDE which is really
a general class of methods for mapping systems $\Sys$ to $\Syshat$ (i.e.\ \TargetClass).
The algorithm ${\mbox{\MapDE}}(\Sys, \Syshat, \Psi)$
returns the system of mapping equations in \rif-form, and, if their integration is successful,
an explicit form of the transformations to map the system to the \TargetClass.
It is described in the {\MapDE} Algorithm ~\ref{Alg:MapDE} provided next.
In that algorithm we suppose that $\hat{\mathcal{L}}, \hat{\mathcal{L}}$ are respectively the Lie algebras of 
symmetries of $\Sys$, $\Syshat$, with defining systems $S$, $\hat{S}$.

For mathematical properties of the algorithms, including finiteness, see the following references.
For {\LAVF} see \cite{LisH:Alg}, for {\rif}'s existence and uniqueness theory see 
\cite{Rus99:Exi}, for the classification of differential rankings see \cite{Rus97:Ran}.
For the algorithmic determination of structure of transitive Lie pseudogroups see Lisle and Reid \cite{LisleReidInfinite98}.

\begin{algorithm}[h]
	\caption{  \MapDE  }  
	\begin{algorithmic}[1]
		\Statex ${\mbox{\MapDE}}(\Source, \Target, \Map)$ 
		\Statex {\bf Input}: 
		\Statex \Source:\, a {\DPS} system $\Sys$, $[x,u]$, $[\xi, \eta]$, Opt
		\Statex \Target:\, a {\DPS} system $\Syshat$, $[\hat{x}, \hat{u}]$, $[\hat{\xi}, \hat{\eta}]$,  Opt
		\Statex \Map:\,   $\Psi$, Opt
		\Statex {\bf Output}:  $\emptyset$ if no consistent \rif-form cases computed, otherwise
		a consistent \rif-form case $Q$ for $\Psi$ and \maple{pdsolve}(Q)		
		\State		Compute															 
		\Statex		\hskip12pt $\Sys :=\mbox{\rif}(\Sys)$, $\mbox{ID}(\Sys)$, $\dim(\Sys)$ 
		\Statex		\hskip12pt $\Syshat :=\mbox{\rif}(\Syshat)$, $\mbox{ID}(\Syshat$), $\dim(\Syshat)$		
		\State		\textbf{if} \maple{DimEquivTest}$(\Sys, \Syshat)\not=$ true \textbf{then} \textbf{return} $\Psi = \emptyset$ \textbf{end if} 
		\State		Compute								 
		\Statex 		\hskip12pt $S = \mbox{\rif}(\mbox{DetSys}(\Sys))$,  $\hat{S} = \mbox{\rif}(\mbox{DetSys}(\Syshat))$
		\Statex		\hskip12pt $\mbox{ID}(S )$, $\mbox{ID}(\hat{S})$, $\dim(S)$, $\dim(\hat{S})$
		\State		\textbf{if} DimEquivTest$(S, \hat{S})\not=$ true \textbf{then} \textbf{return} $\Psi = \emptyset$ \textbf{end if}  
		\State		Compute $\mbox{StrucCons}(\mathcal{L})$, $\mbox{StrucCons}(\hat{\mathcal{L}})$                                
		\Statex		\hskip12pt where $d = \dim(\mathcal{L}) = \dim(\hat{\mathcal{L}})$.
		\State		\textbf{if} $\mathcal{L }\not \simeq \hat{\mathcal{L}}$  \textbf{then}   \textbf{return} $\Psi = \emptyset$ \textbf{end if}      
		\State 		Set $M_{\mbox{\BK}}=$\eqref{eq:BK} and obtain the mapping system:									
		\Statex		\hskip12pt $M := S \cup \; \hat{S}|_{\Psi} \; \cup \; M_{\mbox{\BK}}(\mathcal{L}, \hat{\mathcal{L}}) \; \cup \{\mbox{Det} \mbox{Jac}(\Psi) \not = 0\}  $          
		\Statex		\hskip12pt where $\hat{S}|_{\Psi}$ is $\hat{S}$ evaluated in terms of $(x,u)$ and $(\xi, \eta, \psi, \phi)$ as
		\Statex		\hskip12pt  functions of $(x,u)$ via $\Psi$.
		\State		Compute $M_{\mbox{\rif}} := \text{\rif}(M, \prec, casesplit, mindim = d)$                                                     
		\State		\textbf{if} $M_{\mbox{\rif}}=\emptyset$  \textbf{then return} $\Psi = \emptyset$ \textbf{end if}				
		\State		\textbf{if} $\exists M_{\mbox{\rif}}[\ell] \in  M_{\mbox{\rif}}$ with $d<\infty$ dimensional ID for $\Psi$
		\Statex		\hskip12pt \textbf{then return} $M_{\mbox{\rif}}[\ell]$ and $\maple{pdsolve}(M_{\mbox{\rif}}[\ell])$ \textbf{end if }     
		\Statex		\textbf{end if }               							   
		\State	      Compute $\mbox{Sys}(\Psi) := \textbf{EquivDetSys}(\Sys, \Syshat)$ 										   
		\Statex		 Q := $\emptyset$
		\State		\textbf{while} Q = $\emptyset$  for each consistent sys $ M_{\mbox{\rif}}[k] \in M_{\mbox{\rif}}$    \textbf{do} 
		\Statex		\hskip12pt  $Q := \text{SelSys} (\text{\rif}(\mbox{Sys}(\Psi) \cup M_{\mbox{\rif}}[k], casesplit, mindim = d))$ 
		\Statex		\textbf{end do}		
		\State		\textbf{if} $Q \not = \emptyset$ \textbf{return} $Q$ and \maple{pdsolve}(Q) \textbf{else return} $\Psi = \emptyset$  \textbf{end if }  
	\end{algorithmic}
	\label{Alg:MapDE}
\end{algorithm}


\paragraph{\textbf{Notes for the MapDE Algorithm}}

\begin{itemize}
	\item[Input:] The input \Source $\Sys$ consists differential polynomial system ({\DPS}) of differential polynomials 
	with coefficients in some computable field (e.g.\ $\mathbb{Q}$);
	Opts are additional Options such as input rankings if not default.
	
	\item[Output:]   \maple{pdsolve} is {\Maple} general purpose exact {\PDE} solver: 
	the application of {\Maple}'s \maple{pdsolve} which can not guarantee successful integration of {\DE}.
	
	\item[Step 1:]  Here and throughout {\rif} and {\alg{ID}} refer to {\Maple}'s {\maple{DEtools}} package commands {\maple{rifsimp}} and {\maple{initialdata}} commands.
	Alternatively one could use other Maple packages such as {\maple{diffalg}} or {\maple{DifferentialThomas}}.
	
	\item[Step 2:]  As introduced in \S \ref{sec:EquivDetSysDimEquivTest},  $\maple{DimEquivTest}(\Sys, \Syshat)$ is a simple algorithm for checking some necessary conditions for the existence of a mapping:
	the simplest being $\dim(\Sys) = \dim(\Syshat)$, and include others corresponding to coefficients of the Differential Hilbert Series for ${\Sys}$ and ${\Syshat}$.
	
	\item[Steps 3, 4:]  Restriction to a subalgebra is also possible and can improve efficiency.  Similarly to Step 2, invariant dimension information can lead to early rejection of existence of a mapping:
	the first being that $\dim(S) = \dim(\hat{S})$.
	
	\item[Steps 5, 6:]  {\LAVF} command \maple{StructureCoefficients}
	algorithmically determines the structure constants of the algebras for $d < \infty$.
	{\Maple}'s \maple{LieAlgebras} and
	\maple{DifferentialGeometry} packages, are then used to
	generate the polynomial system for $b_{i,j}$ in a change of basis matrix $B=[b_{i,j}]$
	which is then analyzed by the solver \maple{Triangularize}.
	
	\item[Step 7:]  The change of coordinates to compute $\hat{S}|_{\Psi}$ is accomplished by applying the Maple command \maple{dchange} and using the transformation properties of Lie vector fields \cite{Blu10:App, Olv93:App}.
	
	\item[Step 8:]  Differential elimination with casesplitting is applied and useless computations on branches with ID $<$ mindim$ = d$ wrt
	$(\xi, \eta, \hat{\xi}, \hat{\eta})$ avoided. The ranking $\prec$ ranks the map variables $\Psi$ less than any derivative of the infinitesimals $(\xi, \eta, \hat{\xi}, \hat{\eta})$ yielding
	an uncoupled system in $\Psi$ whose ID is then examined and cases with less than $d$ dimensional data rejected.
	This ranking means that the linearity in $(\xi, \eta, \hat{\xi}, \hat{\eta})$ is maintained in computations.

	\item[Step 11:]  See \S \ref{sec:EquivDetSysDimEquivTest}.
	
	\item[Step 12:]  \maple{SelSys}$(M)$ selects a consistent system from the output of \rif$(M, \text{casesplit}, mindim= d)$

\end{itemize}

\subsection{MapDE for mapping Linear Homogeneous DE to Constant Coefficient Linear DE}
\label{sec:ToConstantCoeffDE}

Here we consider how to map a linear homogeneous source {\DE} to a constant coefficient linear homogeneous {\DE}, 
with an algorithm which results from straightforward changes to Algorithm \ref{Alg:MapDE}.

The idea introduced in Bluman et al.\ \cite[\S 2.5]{Blu10:App} for this problem is to introduce a chain of
Lie subalgebras whose purpose is to focus on the \Target:
$\hat{\mathcal{L}} \supset \hat{\mathcal{L}}^{'}   \supset \hat{\mathcal{L}}^*$
and via $\Psi$ also a chain $\mathcal{L} \supset {\mathcal{L}}^{'}   \supset {\mathcal{L}}^{*}$.

Now $\Syshat = \sum_{i \in I} a_i K(i)= 0$ where $K= \{ K(i) : i \in I \}$ is the set of derivatives of $\hat{u}$ of order $ \leq $ differential order of $\Syshat$.
The unspecified constants $a_i$ are the coefficients of the target.
It is natural to restrict to  transformations that preserve the linearity an homogeneity of the input {\DE} and 
result from eliminating the superposition symmetry:
$\hat{x} = f(x)$ and $\hat{u} = g(x)\,u$ where $x = (x^1, \ldots, x^n)$ (see Bluman et al.\ \cite{Blu10:App}).
Correspondingly its natural to consider a subalgebra $\mathcal{L'}$ that results 
by appending the equations $\xi^j_u = 0,   j = 1, \ldots ,  n$ and $ \eta_u = \eta/u$ to $S$ to form $S'$:
\begin{equation}
\label{eq:SpEq}
S' := \{  \eta_u = \eta/u , \xi^j_u = 0 :   j = 1, \ldots ,  n  \} \cup S
\end{equation}
and similarly for $\hat{\mathcal{L}}^{'}$.
To avoid the early calculations that involve $\Syshat$, we focus like Bluman et al, on accessible infinitesimal  information
encoded in a Lie algebra $\hat{\mathcal{L}}^*$.   In this case $\hat{\mathcal{L}}^*$  corresponds to $n$ commuting translations in the independent variables $\hat{x}^1, \ldots , \hat{x}^n$, i.e.\ 
$n$ translations with generators $\frac{\partial}{\partial \hat{x}^j}$.
The corresponding differential system for $\hat{S}^* $ and $\hat{\mathcal{L}}^*$ is
\begin{equation}
\label{eq:TargetSysConsCoeff}
\hat{S}^* = \{  \hat{\eta} = 0, \hat{\xi}^j_{\hat{u}} = 0,  \hat{\xi}^j_{\hat{x}^k} = 0 : 1 \leq j,\, k \leq n \} 
\end{equation}





\begin{algorithm}[h]
	\caption{  \MapDE with \maple{TargetClass}=\maple{ConstantCoeffDE}}  
	\begin{algorithmic}[1]
		\Statex ${\mbox{\MapDE}}(\Source, \Target, \Map)$ 
		\Statex {\bf Input}: 
		\Statex \Source:\, A single linear homogeneous {\DPS} $\Sys$, $[x,u]$, $[\xi, \eta]$,
		\Statex \hskip12pt Opt
		\Statex \Target:\, \maple{TargetClass}=\maple{ConstantCoeffDE} , Opt
		\Statex \Map:\,   $\Psi$, Opt
		\Statex {\bf Output}: $\emptyset$ if no consistent \rif-form cases computed, otherwise
		a consistent \rif-form case $M_{\mbox{\rif}}[\ell]$ for $\Psi$ and  $\maple{pdsolve}(M_{\mbox{\rif}}[\ell])$
		\State		Compute															 
		\Statex		\hskip12pt $\Sys :=\mbox{\rif}(\Sys)$, $\mbox{ID}(\Sys)$, $\dim(\Sys)$  
		\Statex \hskip12pt $\Syshat := \sum_{i \in I} a_i K(i)= 0$
		
		\State 	 Compute $S = \mbox{\rif}(S')$,  $\mbox{ID}(S' )$, $\dim(S)$, $\dim(\hat{S})$	
		\State Compute $\mbox{StrucCons}(\mathcal{L}')$
		
		\State $M := S' \cup \; \hat{S}^*|_{\Psi} \; \cup \; M_{\mbox{\BK}}(\mathcal{L}', \hat{\mathcal{L}}^* ) \; \cup \{\mbox{Det} \mbox{Jac}(\Psi) \not = 0\}  $ 
		\State 
		Compute $M_{\mbox{\rif}} := \text{\rif}(M, \prec, casesplit, mindim = n)$                                       
		\State		\textbf{if} $M_{\mbox{\rif}}=\emptyset$  \textbf{then return} $\Psi = \emptyset$ \textbf{end if}	                                    
		\State		\textbf{if} $\exists M_{\mbox{\rif}}[\ell] \in  M_{\mbox{\rif}}$ with $d=n<\infty$ dimensional ID for $\Psi$ 
		\Statex	\hskip12pt \textbf{then return} $M_{\mbox{\rif}}[\ell]$ and $\maple{pdsolve}(M_{\mbox{\rif}}[\ell])$
		\Statex	\textbf{elif return} $\Psi = \emptyset$
		\Statex	\textbf{end if }    
		
	\end{algorithmic}
	\label{Alg:MapDECoeffOption}
\end{algorithm}

\section{Examples}
\label{sec:Examples}
In this section we apply our algorithm to examples. 

\subsection{Equivalence}
\label{se:Equivalence}

\begin{example}
	Bluman et al.\ \cite[\S 2.3.2, pg 133-137]{Blu10:App} apply their mapping method based on explicit integrations to determine an invertible mapping by a point transformation of the cylindrical KdV equation $\Sys$ to the KdV equation $\Syshat$ that first appeared in the work of Korobeinikov \cite{Korobeinikov82}:
	\begin{align}
		\label{eq:KdVsys}
		\Sys      &:= \left\{u_{{x,x,x}} =-uu_{{x}}-u_{{t}}-{\frac {u}{2t}} \right\}  \\
		\Syshat &:=  \left\{ {\it \uh}_{{{\it \xh},{\it \xh},{\it \xh}}} =-{\it \uh}\,{\it \uh}_{{{\it \xh}}}-{\it \uh}_{{{\it \thn}}} \right\}
	\end{align}
	They give details of their calculations and for illustration we apply the \MapDE algorithm \ref{Alg:MapDE} to the same example.
	Here we seek transformations $\hat{x} = \psi(x,t,u),\hat{t} = \phi(x,t,u), \hat{u} = \varUpsilon(x,t,u)$.
	
	\begin{itemize}
		\item[Steps 1, 2:] Both $\Sys$ and $\Syshat$ are already in \rif-form with respect to any orderly ranking.
		The initial data for the $\Sys$ and $\Syshat$ are
		\begin{align*}
			\{u(x_{0},t)={\it F_1} \left( t \right) , u_{x}(x_{0},t)={\it F_2} \left( t
			\right) ,  u_{x,x}(x_{0},t)={\it F_3} \left( t \right) \}
		\end{align*}
		\begin{align*}
			\{
			{\it \uh}(\xh_{0}, \thn)={\it F_1} \left( {\it \thn} \right), {\it \uh}_{\it \xh}(\xh_{0},\thn)={\it F_2} \left( {\it \thn} \right) , {\it \uh}_{\xh, \xh}(\xh_{0}, \thn)={\it F_3} \left( {\it \thn} \right) \}
		\end{align*}
		Here there are arbitrary functions in the initial data, \\so $\dim{\Sys} = \dim{\Syshat} = \infty$.  
		Their Hilbert Series obviously are equal, and up to the order of involutivity: $\mbox{HS}(s) = 1 + 2s + 3s^2 + 0(s^3)$ where the coefficient of $s^n$ is the number of parametric derivatives of order $n$.
		So $\mbox{DimEquivTest}(\Sys, \Syshat) = \mbox{true}$ in Step 2.
		
		\item[Step 3, 4:] The \rif-form systems $S, \hat{S}$ are 
		\begin{align}
			S =[\eta_{{u,u}}&=0, \xi_{{x}}=-\frac{1}{2}\,\eta_{{u}}, \beta_{{x}}=0,  \eta_{{x}}={\frac{-3\,t\eta_{{u}}-2\,\beta}{4{t}^{2}}},  \qquad \quad  \nonumber  \\  
			\beta_{{t}}&=-\frac{3}{2}\,\eta_{{u}}, \eta_{{t}}={\frac {4\,t\eta_{{u}}u+2\,\beta\,u-\eta\,t}{2{t}^{2}}},   \nonumber \\
			\xi_{{u}}&=0, \beta_{{u}}=0,   \xi_{{t}}=-\eta_{{u}}u+\eta ]  
		\end{align} 
		\begin{align}	
			\hat{S} = [{ \hat{\eta}}_{{{ \hat{u}},{ \hat{u}}}}&=0,{ \hat{\xi}}_{{{ \hat{x}}}}=-\frac{1}{2}{\hat{\eta}}_{{{ \hat{u}}}},{ \hat{\beta}}_{{{ \hat{x}}}}=0,{ \hat{\eta}}_{{{ \hat{x}}}}=0 ,   \\
			{\hat{\xi}}_{{{ \hat{t}}}}&=-{ \hat{\eta}}_{{{ \hat{u}}}}{ \hat{u}}+{ \hat{\eta}},{ \hat{\beta}}_{{{\hat{t}}}}=-3/2\,{ \hat{\eta}}_{{{ \hat{u}}}},{ \hat{\eta}}_{{{\hat{t}}}}=0,{ \hat{\xi}}_{{{ \hat{u}}}}=0,{ \hat{\beta}}_{{{ \hat{u}}}}=0] \nonumber
		\end{align}
		and yield $\mbox{ID}$ giving $\dim(S) = \dim(\hat{S}) = 5$.  Also DimEquivTest$(S, \hat{S})=$ true in Step 4.
		\item[Step 5:]  Here {\MapDE} uses the {\LAVF} command \maple{StructureConstants}  to compute the structure of the $4$ dimensional Lie algebras for $\Sys$ and $\Syshat$ obtaining:
		\begin{align}
			\label{eq:StructKdv}
			\mathcal{L}:   &  [  [Y_{{1}},Y_{{4}}]=-1/2\,Y_{{1}},[Y_{{2}},Y_{{3}}]=Y_{{1}},[Y_{{2}},Y_
			{{4}}]=-3/2\,Y_{{1}}+Y_{{2}},    \nonumber  \\
			& [Y_{{3}},Y_{{4}}]=-3/2\,Y_{{3}},[Y_{{1}}, Y_{{3}}]=0,[Y_{{1}},Y_{{2}}]=0  ]  \\
			\hat{\mathcal{L}}:  & [ [{\it \hat{Y}}_{{1}},{\it \hat{Y}}_{{4}}]=-1/2\,{\it \hat{Y}}_{{1}},[{\it \hat{Y}}_{{2}},
			{\it \hat{Y}}_{{3}}]={\it \hat{Y}}_{{1}},[{\it \hat{Y}}_{{2}},{\it \hat{Y}}_{{4}}]=-3/2\,{
				\it \hat{Y}}_{{2}},  \nonumber  \\
			&[{\it \hat{Y}}_{{3}},{\it \hat{Y}}_{{4}}]={\it \hat{Y}}_{{3}},[{\it \hat{Y}}
			_{{1}},{\it \hat{Y}}_{{2}}]=0,[{\it \hat{Y}}_{{1}},{\it \hat{Y}}_{{3}}]=0  ]
		\end{align}
		
		\item[Steps 5, 6:]  We obtain $\mathcal{L} \simeq \hat{\mathcal{L}}$ and the explicit isomorphism:
		\begin{equation}
		{\it \hat{Y}}_{{1}}=Y_{{1}},\quad {\it \hat{Y}}_{{2}}=Y_{{3}},\quad {\it \hat{Y}}_{{3}}=Y_{{1}}-
		Y_{{2}},\quad {\it \hat{Y}}_{{4}}=Y_{{4}}
		\end{equation}
		Bluman et al.\ \cite[Eqs (2.39), (2.40), pg 134]{Blu10:App}obtain the structure and an isomorphism by explicitly integrating the defining systems,
		whereas we avoid this.
		This isomorphism is a necessary but not sufficient condition for the existence of a local analytic
		invertible map to $\Syshat$.
		
		\item[Steps 7, 8:]  The \rif-form of the mapping system results in one consistent case with $d=4$-dim ID in
		the infinitesimals for $S, \hat{S}$.  The \rif-form of the $\Psi$ system is:
		\begin{align}
			\label{eq:PsiKdv}
			[\phi_{{t,t}} &=-3/2\,{\frac {\phi_{{t}}}{t}},\varUpsilon_{{u,u}}=0,
			\varUpsilon_{{x}}=-1/2\,{\frac {\varUpsilon_{{u}}}{t}},\psi_{{x}}=\phi_{{t}} \varUpsilon_{{u}},  \nonumber \\ \phi_{{x}}&=0, 
			\varUpsilon_{{t}}={\frac {\varUpsilon_{{u}}u}{t}
			},\psi_{{t}}=-u\varUpsilon_{{u}}\phi_{{t}}+\varUpsilon\,\phi_{{t}},\psi_{{u}
			}=0,\phi_{{u}}=0]
		\end{align}
		
		\item[Steps 9, 10:]  Step 9 does not apply. The above \rif-form for the $\Psi$ system has ID for $z_0= (x_0, t_0, u_0)$:
		$$  \varUpsilon(z_0) = c_1, \varUpsilon_u(z_0) = c_2, \phi(z_0) = c_3,  \phi_t(z_0) = c_4, \psi(z_0) = c_5 $$
		Geometrically, since the ID has dimension $5 > d=4$, there is class of systems with the same $d$ dimensional 
		invariance group, that possibly includes the \Target system.  So Step 10 does not apply. Thus we have to apply 
		$\maple{EquivDetSys}$ to find missing condition(s).  After explicit integration Bluman et al also find that they don't 
		uniquely specify the target, and essentially they substitute the transformations to obtain the parameter values to 
		specify the target.  
		
		\item[Steps 11, 12:]  Applying \rif to the combined system  $\{ \textbf{EquivDetSys}(\Sys, \Syshat), (\ref{eq:PsiKdv})\}$ 
		yields a single case:
		\begin{align*}
			M_{\mbox{\rif}}= [\psi_{{t,t}}&=-3/2\,{\frac {\psi_{{t}}}{t}},\quad \varUpsilon_{{u,u}}=0, \quad
			\varUpsilon_{{x}}=-1/2\,{\frac {\varUpsilon_{{u}}}{t}},\quad \\ \phi_{{x}}&=0,\quad 
			\psi_{{x}}={\frac{\varUpsilon_{{u}}\psi_{{t}}}{-u\varUpsilon_{{u}}+\varUpsilon}}, \quad
			\varUpsilon_{{t}}={\frac {u\varUpsilon_{{u}}}{t}}, \quad \\ \phi_{{t}}&={\frac {\psi_{{t}}}{-u\varUpsilon_{{u}}+\varUpsilon}},\quad \phi_{{u}}=0,\quad \psi_{{u}}=0]
		\end{align*}
		where the constraint is $
		-{\psi_{{t}}}^{2}{\varUpsilon_{{u}}}^{3}+{\varUpsilon_{{u}}}^{2}{u}^{2}-2\,\varUpsilon\,u\varUpsilon_{{u}}+{\varUpsilon}^{2}=0
		$,
		and the inequation ${\frac {{\psi_{{t}}}^{2}{\varUpsilon_{{u}}}^{2}}{ \left( -u\varUpsilon_{{u}}+\varUpsilon \right) ^{2}}} \not = 0$.
		The ID shows we now have $4 = d$ parameters, confirming the existence of the transformations, without integration.

		\item[Step 13]
		Applying \maple{pdsolve} yields the solution for the transformation below:
		
		\begin{align*}
			\{{\it \xh}={\frac {{\it c_3}\,x}{\sqrt {t}}}-2\,{\frac {{\it c_3}\,{\it c_2}}{\sqrt {t}{\it c_1}}}+{\it c_4}, \quad
			{\it \thn}=-2\,{\frac {{\it c_3}}{\sqrt {t}{\it c_1}}}+{\it c_5}, \quad \\ {\it \uh}  =1/2\, \left( 2\,tu-x \right) {\it c_1}+{\it c_2}\}
		\end{align*}
		where
		\begin{align*}
			\left( -1/4\,{{\it c_1}}^{3}{{\it c_3}}^{2}+1/4\,{{\it c_1}}^{2} \right) {x}^{2}+ \left( {{\it c_1}}^{2}{\it c_2}\,{{\it c_3}}^{2}-{\it c_1}\,{\it c_2} \right) x\\
			-{\it c_1}\,{{\it c_2}}^{2}{{\it c_3}}^{2}+{{\it c_2}}^{2}=0.
		\end{align*}
		subject to the determinental condition.  The last condition implies $c_1 c_3^2 = 1$.  Specializing the values 
		of the $c_j$ give the transformations obtained also in Bluman et al.
		
	\end{itemize}
	
\end{example}

\subsection{Mapping to constant coefficient DE}
\label{sec:PDEEx}

\begin{example}
	The harmonic-oscillator Schr\"odinger Equation $i \hslash  \varphi_t = -\frac{\hslash^2}{2 m} \varphi_{xx} +\frac{1}{2}m \omega^2 x^2 \varphi $
	which in normalized \rif-form is:
	\begin{equation}
	\label{ex:ConsCoeffR}
	\Sys := \{ u_{{x,x}}=-{x}^{2} u +u_{{t}} \} 
	\end{equation}
	Applying the $\MapDE$ algorithm using the option \maple{Target = ConstantCoeffDE} shows that  (\ref{ex:ConsCoeffR}) maps to a constant coefficient linear {\DE}: 
	\begin{equation}
	\label{ex:ConsCoeffRh}
	\Syshat:= \left\{ a_1 {\it \uh}_{{{\it \xh},{\it \xh}}}+ a_2 \,{\it \uh}_{{{\it \xh},{\it \thn}}}+a_3 {\it \uh}_{{{\it \thn},{\it \thn}}}  + a_4  {\it \uh}_{{{\it \xh}}}  + a_5    {\it \uh}_{{{\it \thn}}} + a_6 {\it \uh}     =0 \right\} 
	\end{equation}
	Existence of such a mapping is given algorithmically and the output includes the system for $\Psi$ with $\dim(\Psi) = 6$:
	\begin{align*}
		\label{eq:conscoeffPsi}
		\Psi = [\psi_{{u}}&=0,\varUpsilon_{{u}}=0,\varUpsilon_{{x,x}}=0,\phi_{{u}}={\frac {\phi}{u}},\psi_{{x}}=\varUpsilon_{{x}},\psi_{{t}}=-{\varUpsilon_{{x}}}^{2}+ \varUpsilon_{{t}},   \nonumber    \\ 
		\phi_{{x}}&=-\,{\frac {\phi\, \left( -{\varUpsilon_{{x}}}^{2}+\varUpsilon_{{t}} \right) }{2 \varUpsilon_{{x}}}}, 
		\varUpsilon_{{t,t,t}}= -\,{\frac {12\,\varUpsilon_{{t,t}}x\varUpsilon_{{x}}-16\,{\varUpsilon_{{t}}}^{2}-3\,{\varUpsilon_{{t,t}}}^{2}}{2\varUpsilon_{{t}}}},   \nonumber    \\
		\varUpsilon_{{x,t}}&= -\,{\frac {\varUpsilon_{{x}} \left( 4\,x\varUpsilon_{{x}}-\varUpsilon_{{t,t}} \right) }{2 \varUpsilon_{{t}}}},       \\ 
		\phi_{{t}}&={\frac {\phi\, \left(4\,{\varUpsilon_{{x}}}^{3}x - 4\,{\varUpsilon_{{x}}}^{2}\varUpsilon_{{t}}{x}^{2}-{\varUpsilon_{{x}}}^{4}\varUpsilon_{{t}}+2\,{\varUpsilon_{{x}}}^{2}{\varUpsilon_{{t}}}^{2}-{\varUpsilon_{{t}}}^{3}-{\varUpsilon_{{x}}}^{2}\varUpsilon_
				{{t,t}} \right) }{4{\varUpsilon_{{x}}}^{2}\varUpsilon_{{t}}}}]   \nonumber 
	\end{align*}
	Integrating the system and specializing the $6$ constants gives:
	\begin{align*}
		{\it \thn} &= \psi(x,t,u) =\frac{2(x+\cos(2t))}{\sin(2t)}, \quad
		{\it \xh}  = \varUpsilon(x,t,u) =\frac{2x}{\sin(2t) },         \\
		{\it \uh} &= \phi(x,t,u) =  \sqrt{\sin(2t)} \exp \left[ \frac{ 2(x^2+1) \cos(t)^2-(x-1)^2}{2 \sin(2t)}  \right] \cdot  u    
	\end{align*}
	where $\Syshat$ is  
	$  {\it \uh}_{{{\it \xh},{\it \xh}}}+ 2 \,{\it \uh}_{{{\it \xh},{\it \thn}}}+ {\it \uh}_{{{\it \thn},{\it \thn}}}  
	-   {\it \uh}_{{{\it \thn}}}     =  0$.
	

	When the mapping system
	$M := S' \cup \; \hat{S}^*|_{\Psi} \; \cup \; M_{\mbox{\BK}}(\mathcal{L}', \hat{\mathcal{L}}^* ) \; \cup \{\mbox{Det} \mbox{Jac}(\Psi) \not = 0\}  $  
	is reduced to \rif-form as in Step 4, of Algorithm \ref{Alg:MapDECoeffOption} it yields a consistent system for $\Psi$ with $10$ arbitrary constants in its initial data, hence establishing
	existence of a mapping to a constant coefficient {\DE}.  
	Since $\dim(\mathcal{L}) = 6$ this means that there is a $4$ dimensional target class of constant coefficient
	linear {\DE}.  Indeed there are two options here, one would be attempt to (heuristically) integrate this system.  We did try this, and succeeded to find the solutions with 10 parameters.  Instead in the spirit of our approach in this paper, to reduce as much as possible to algorithmic differential elimination, we also experimented with 
	inserting a step in the algorithm, that executes a general change of coordinates from the \Target {\DE} which has $5$ arbitrary constants, ranking those constants highest in the ordering and obtained as a subsystem:

	\begin{align*}
		a_{{2}}&=2\,{\frac {\varUpsilon_{{x}}}{\psi_{{x}}}}, \, 
		a_{{3}}={\frac{{\varUpsilon_{{x}}}^{2}}{{\psi_{{x}}}^{2}}}, \,
		a_{{4}}={\frac{-\phi\,\psi_{{t}}-2\,\psi_{{x}}\phi_{{x}}}{{\psi_{{x}}}^{2}\phi}}, \, 
		a_{{5}}=-{\frac {\phi\,\varUpsilon_{{t}}+2\,\varUpsilon_{{x}}\phi_{{x}}}{{\psi_{{x}}}^{2}\phi}},\\  
		a_{{6}} &=1/4\,{\frac {4\,{\phi}^{2}{x}^{2}\varUpsilon_{{t}}-4\,{\phi}^{2}x\varUpsilon_{{x}}+{\phi}^{2}\varUpsilon_{{t,t}}+4\,\phi\,\varUpsilon_{{t}}\phi_{{t}}+4\,\varUpsilon_{{t}}{\phi_{{x}}}^{2}}{\varUpsilon_{{t}}{\psi_{{x}}}^{2}{\phi}^{2}}} 
	\end{align*}
	The transforming system $M_{\mbox{\rif}}$ is updated by adding these equations. 
	So the constants amount to integrations, and can be used to reduce the dimension of the system for $\Psi$ by specializing their values.
	This yielded the final $6$ dimensional system.  Note that it was reduced by $4$ dimensions (after elimination the relation 
	$a_2^2  = 4 a_3 $).
	
\end{example}

\begin{example}
	\label{ex:ToConstantCoeffLinear}
	Consider the DE arises in financial models known as Black Schole's\cite{MasMa119:ExBlacSch} as the \Source system $\Sys$,
	\[
	v_{{t}} + \frac{s^{2}v_{s,s}}{2}+sv_{{s}}-v=0
	\]
	Using our algorithm $\MapDE$ with \maple{TargetClass = ConstantCoeffDE}, automatically yields the mapping $\Psi$:
	\[
	{\it {\shatn}}=\ln  \left( s \right) ,\, {\it {\thatn}}= t, \, {\it {\vhatn}}=v
	\]	
	and the \Target system $\Syshat$ is ${\it {\vhatn}}_{{{\it {\thatn}}}}+1/2\,{\it {\vhatn}}_{{{\it {\shatn}},{\it {\shatn}}
	}}-{\it {\vhatn}}=0$.
	
\end{example}
\section{Discussion}
\label{sec:Discussion}

Mappings of mathematical models are a fundamental tool of mathematics and its applications.
This fact and the notorious difficulty of their computation motivates us to explore the approach we presented in this article.
This resulted in our algorithm, \MapDE, which is given algorithmic realization for two cases.
The first is where the input systems $\Sys$ and $\Syshat$ are specified as polynomially nonlinear DEs and \MapDE returns a reduced involutive \rif-form for the mapping equations.
The second is where $\Sys$ is a linear homogeneous {\DE} and the \TargetClass is a constant coefficient linear 
homogeneous {\DE} (\maple{Target = ConstantCoeffDE}).  A key aspect of our approach
is to exploit the linearity of the Bluman-Kumei mapping equations that arise in the presence of symmetry, and postpone, 
simplify and even avoid direct computations with the full nonlinear determining equations for the mappings.
We also implement some fast preliminary tests for equivalence under mappings.

The closest approach to the our work, are the works of Bluman and collaborators and in particular the work by 
Anco, Bluman and Wolf \cite{AnBlWo110:Mapp} and also Wolf \cite{Wolf116:SymSof}, which considered a computer program for computing linearization mappings.
It exploits Wolf's program {\sc{ConLaw}}'s strong facilities for integrating systems of {\PDE} exactly in addition to the {\BK} mapping equations, as well as an embedding technique involving multipliers and conservation laws.
They also mention that the problem of full algorithmization using differential algebra as an important open problem.


In another paper \cite{MohamadiReidHuang19}, we give various extensions {\MapDE}.
One of these involves extending it to determining existence of exact linearization mappings of {\DE}.  In so doing we provide algorithm and combining aspects of the approach of Bluman, Anco and Wolf \cite{AnBlWo110:Mapp,Wolf116:SymSof} and also of Gerdt et al \cite{LGM101:LG}.

Building in invariant properties into the completion process is also a possibility; borrowing aspects of the more 
geometrical approaches.

Longer term we are particularly interested in exploring approximate mappings and approximate equivalence.
Indeed {\LAVF} already has the first available algorithm for determining the structure of approximate symmetry
of {\DE} theoretically first described in Lisle, Huang and Reid \cite{LisHG:Sym}.  Indeed consider Poisson's equation for a gravitational 
potential $u(x,y,z)$: $\nabla^2 u = f(x,y,z)$ and an interstellar gas with density proportional 
to $f(x,y,z) = \frac{1}{2}(G(x,y,z-a) + G(x,y,z+a)) $ where $G(x,y,z) = \exp(-x^2-y^2-z^2)$ and $a=10^{-3}$:
$$u_{xx} + u_{yy} + u_{zz} = f(x,y,z) = \frac{1}{2} \left( G(x, y, z+a) + G(x,y,z-a) \right)$$

Applying 
Lie's standard method where 
$$L = \xi(x,y,z,u)\dd{x} +  \eta(x,y,z,u)\dd{y} + \psi(x,y,z,u)\dd{z} + \phi(x,y,z,u)\dd{u}$$
Discarding the superposition via $\phi_u = \phi/u$ and performing an exact
symmetry analysis yields only a $1$ dimensional rotation group about the $z$ axis;
throughout all space no matter how small $a$ is.

We now apply our method to find approximate Lie algebra of symmetries of Poisson's Equation
at the point $(x,y,z) = (0,3.2,0)$ and $a = 10^{-3}$.
Recall we only found a $1$ dimensional Lie symmetry group of rotations in the $x-y$ plane about the $z$ axis.
We find:
\begin{eqnarray*}
	[L_{{1}}, L_{{2}} ] &=& -{ 1.182\times 10^{-13}}\, L_{{1}}-  { 2.724\times 10^{-9}}\, L_{{2}}- 0.707 \, L_{{3}}    \\    
	\phantom{YY}  [  L_{{1}}, L_{{3}} ]  &=&  - { 1.446\times 10^{-7}}\, L_{{1}}+ 0.236 \, L_{{2}}+ { 2.724\times 10^{-9}}\, L_{{3}}  \\
	\phantom{YY} [  L_{{2}},  L_{{3}}   ]   &=& - 0.707\, L_{{1}}+{ 1.446\times 10^{-7}} \, L_{{2}}-{ 6.042\times 10^{-14}}\,L_{{3}}   \\
\end{eqnarray*}
which we can recognize as
$$ [L_{{1}},L_{{2}}]=-\frac{1}{\sqrt {2}} L_{{3}}, \hskip6pt [L_{{1}},L_{{3}}]=\frac{1}{3\sqrt {2}} L_{{2}},  \hskip6pt [L_{{2}},L_{{3}}]=-\frac{1}{\sqrt {2}}L_{{1}}  $$
or after the basis change $L_1 = \frac{Y_1}{\sqrt{6}}, L_2 = - \frac{Y_2}{\sqrt{2}}  ,  L_3 = \frac{Y_3}{\sqrt{6}}$ is $so(3)$:
$$ [Y_1 , Y_2 ] = Y_3,  \quad   [Y_2 , Y_3 ] = Y_1,  \quad   [Y_3 , Y_1 ] = Y_2  $$

Remarkably we get different regions with different approximate groups, plus transition bands.  Potentially and intuitively the model can be mapped to various forms depending on the region, which
is a longer term topic, for our research.

\begin{figure}[h!]
	\begin{centering}
		\includegraphics[angle=0,origin=c, width=4cm]{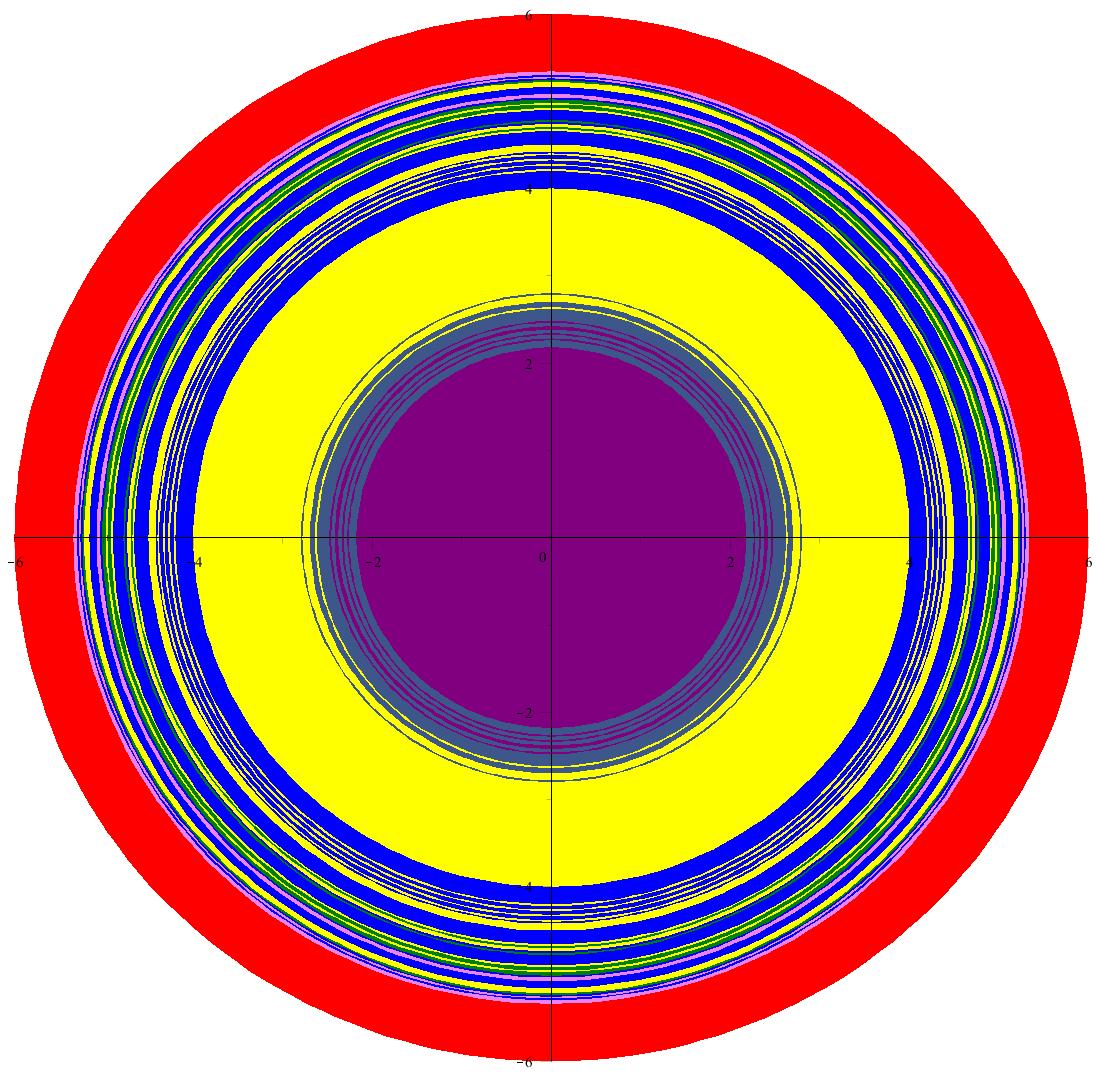}
		\label{fig:Poissonsym}
		\caption{Regions of approximate symmetry in the $x-y$ plane for  $\nabla^2u= f(x,y,z)$.
			Purple Region: $\dim{\mathcal{L}}=1$, ${\mathcal{L}}\approx$so(2); Yellow Region $\dim{\mathcal{L}}=3$,$ {\mathcal{L}}\approx $so(3);  Red Region $\dim{\mathcal{L}}=11$}
	\end{centering}
\end{figure}

\vspace{1cm}
\section*{Acknowledgments}
GJR acknowledges the support of an NSERC Discovery Grant from the government of Canada.
GJR acknowledge the contribution of his colleague and departed friend Ian Lisle, whose inspiration lies behind this work in its spirit and many details.

\end{document}